# 1D+4D-VAR data assimilation of lightning with WRFDA system using nonlinear observation operators

Răzvan Ştefănescu[1], I. Michael Navon[1], Henry Fuelberg[2], Max Marchand[2]

[1]The Florida State University, Department of Scientific Computing, Tallahassee, Florida 32306, USA, rstefanescu@fsu.edu, inavon@fsu.edu

[2]The Florida State University, Department of Earth, Ocean, and Atmospheric Science, Tallahassee, Florida 32306, USA
hfuelberg@fsu.edu, mmarchand@fsu.edu

June 2013

**Abstract**

This paper addresses the impact of assimilating data from the Earth Networks Total Lightning Network (ENTLN) during two cases of severe weather. Data from the ENTLN serve as a substitute for those from the upcoming launch of the GOES Lightning Mapper (GLM). We use the Weather Research and Forecast (WRF) model and variational data assimilation techniques at 9 km spatial resolution. The main goal is to examine the potential impact of lightning observations from the future GLM.

Previous efforts to assimilate lightning observations mainly utilized nudging approaches. We develop three more sophisticated approaches, 3D-VAR WRFDA and 1D+nD-VAR (n=3,4) WRFDA schemes that currently are being considered for operational implementation by the National Centers for Environmental Prediction (NCEP) and the Naval Research Laboratory (NRL). This research uses Convective Available Potential Energy (CAPE) as a proxy between lightning data and model variables.

To test the performance of the aforementioned schemes, we assess the quality of resulting analysis and forecasts of precipitation compared to those from a control experiment and verify them against NCEP stage IV precipitation. Results demonstrate that assimilating lightning observations improves precipitation statistics during the assimilation window and for 3-7 h thereafter. The 1D+4D-VAR approach performs best, significantly improving precipitation root mean square errors by 25% and 27.5% compared to the control during the assimilation window on the two cases. This finding confirms that the variational nD-VAR (n=3,4) assimilation of lightning observations improves the accuracy of moisture analyses and forecasts. Finally, we briefly discuss limitations inherent in the current lightning assimilation schemes, their implications, and possible ways to improve them.

# 1. Introduction

Meteorological satellite data help describe the atmospheric condition, the climate, ocean, and earth system in general by enriching the number and distribution of observations. Current satellite instruments produce large amounts of information that are integrated into data assimilation systems operating in numerical weather prediction (NWP) centers around the world to increase analysis skills and forecast scores.

Advances in science and technology continue to improve the temporal, spatial, and spectral resolution of satellite observations and allow new products to be generated. NASA will launch the GOES-R Lightning Mapper (GLM, Goodman et al. 2013) in 2015 that will provide continuous, full disk, high resolution total lightning data (intracloud (IC) and cloud-to-ground (CG)). Our research explores the capabilities of current state of the art data assimilation techniques to prepare for optimal use of the upcoming GOES-R lightning products.

Lightning data are a relatively new source of information for use in data assimilation. The earliest effort to assimilate lightning data can be attributed to Alexander et al. (1999) who improved simulated latent heating rates using data derived from spaceborne and lightning-derived rainfall measurements to more accurately forecast precipitation associated with an extratropical cyclone. Chang et al. (2001) used continuous lightning distributions tuned with intermittent satellite microwave data through a probability matching technique to continuously generate convective rainfall rates that were assimilated into a mesoscale numerical model. Their simulated moisture characteristics of the 1998 Groundhog Day Storm were more accurate than those without assimilation. Papadopoulos et al. (2005, 2009) used CG lightning information from a European regional network of ground receivers to provide information about the spatiotemporal development and intensity of deep moist convection to a mesoscale meteorological model. They then nudged model generated humidity profiles to empirical profiles as a function of the observed lightning intensity. The assimilation method proposed by Mansell et al. (2007) used lightning data from the Oklahoma Lightning Mapping Array (LMA; Rison et al.1999, Thomas et al. 2004) and the National Lightning Detection Network (NLDN; Biagi et al. 2007, Orville 2008, Cummins and Murphy 2009) as a proxy for the presence or absence of deep convection to control the Kain-Fritsch (KF; Kain and Fritsch 1993) convective parameterization scheme. The World-Wide Lightning Location Network (WWLLN) is a valuable source for detecting both CG and IC flashes which has recently increased its detection efficiency (Rodger et al. 2009).

NOAA ESRL/GSD conducted the first research to utilize lightning data in preliminary operational applications. Ground strike densities were used to supplement radar-derived reflectivity (Weygandt et al. 2008) from the National Severe Storms Laboratory's (NSSL) national (CONUS) 3D radar mosaic which were then assimilated into the Rapid Update Cycle (RUC) model (Benjamin et al. 2004a, b, 2006, 2007). Lightning strike locations were employed to identify convective grid points in Fita et al. (2009). Pessi and Businger (2009a, b) developed a lightning - convective rainfall rate relationship for use in the MM5 model based on long range lightning data and rainfall data from the Tropical Rainfall Measuring Mission (TRMM) satellite. The Lightning Potential Index was used to approximate flash rates in two heavy rain events using the WRF model (Lynn and Yair 2010). Fierro et al. (2012) employed lightning data from the Earth Networks Total Lightning Network (ENTLN) to force deep convection in the WRF model at the locations of observed storms. Results exhibited an improved representation of observed severe convection at the analysis time and during the forecast period.

Most prior efforts to assimilate lightning data into numerical models have focused on Newtonian nudging techniques in which vertical profiles of humidity and/or latent heat were altered. Most of the experiments were performed at horizontal grid spacings greater than 10 km, except for Fierro et al. (2012) whose resolution was at cloud resolving scales.

The present paper uses ENTLN ground based lightning data to investigate the effects of assimilating lightning data at 9 km resolution during two cases of severe storms: a supercell event occurring predominantly in Mississippi and Alabama on 27 April 2011, and a squall line event that initiated in Kentucky and Tennessee and later spread to coastal South Carolina and Georgia on 15 June 2011.

The ENTLN lightning information was clustered over $9 \times 9$ km grid cells to mimic the resolution of the GLM. Thus, they are considered pseudo GLM data. The GLM will map total lightning at a resolution varying between $8 - 12$ km and is expected to achieve a detection efficiency in excess of 70% (e.g., Goodman et al. 2013).

The global and regional NWP data assimilation systems use variational techniques such as 3D-VAR (Brousseau et al. 2008), incremental 3D-VAR (Barker et al. 2004, Derber and Wu 1998, Kleist et al. 2009), incremental 4D-VAR (Rabier et al. 2000, Honda et al. 2005, Gauthier et al. 2007, Laroche et al. 2007, JMA 2007, Rawlins et al. 2007, Huang et al. 2009, Fillion et al. 2010) and hybrid ensemble 3D-VAR formulation (Buehner 2005, Wang et al. 2008a, b). ERA-Interim is the latest global atmospheric reanalysis produced by the European Centre for Medium-Range Weather Forecasts (ECMWF) and it is described in Dee et al. (2011). For assimilation of cloud and precipitation observations approaches like ensemble data assimilation (Zupanski et al. 2011a,b, Zhang et al. 2013), 4D-VAR (Lopez 2011, Benedetti and Janisková 2008, McNally 2009, Karbou et al. 2010, Amerault et al. 2009) and the two steps 1D+4D-VAR (Marecal and Mahfouf 2002, 2003, Bauer et al. 2006a, b, Lopez and Bauer 2007, Geer et al. 2008, Benedetti et al. 2005, Janisková et al. 2012, Wang et al. 2012, Weng et al. 2007) method have been employed. A new radiance-assimilation scheme for microwave-imager observations, which unifies the treatment of clear-sky, cloudy and precipitation-affected situations (Bauer et al. 2010, Geer et al. 2010, Geer and Bauer 2011) is operational in the ECMWF 4D-Var system. A discussion about assimilation of cloud and precipitation data issues can be found in Errico et al. (2007). Bauer et al. (2011) presents an overview of satellite cloud and precipitation assimilation techniques at operational NWP centers.

The novelty of our research consists in applying three complex nD-VAR (n=3,4) variational data assimilation techniques to ingest the pseudo GLM flash rates into the WRF model. The goal is to better represent convection at the analysis time and thereby produce more accurate short term forecasts. The highly nonlinear operator that we employ uses convective available potential energy (CAPE) as a proxy between lightning data and WRF model variables in the WRFDA system. Thus, we exploit power law relationships between total lightning flash rate and maximum vertical velocity as described in Barthe et al. (2010), Deierling et al. (2005, 2008), and Deierling and Petersen (2008). These schemes have been widely used to study the production of nitrogen oxides from lightning (Pickering et al.1998, Fehr et al.2004, Salzmann et al. 2008, Barthe and Barth 2008). Correlations between the occurrence of total lightning activity and deep convective updrafts, described by Goodman et al. (1998), Wiens et al. (2005), and Fierro et al. (2006), also allow the development of lightning proxies within the NWP models.

Our first scheme uses 3D-VAR to directly assimilate lightning observation into the WRF model by adjusting the temperature lapse rate. Marecal and Mahfouf (2003) showed that nonlinearities can cause the direct assimilation of rainfall rates to be unsuccessful in the ECMWF 4D-VAR (using the incremental formulation proposed by Courtier et al. (1994)). In our case we empirically prove that the direct assimilation of lightning data into the WRF 3D-VAR schemes is threshold limited, due to the highly nonlinear observation operator and nature of the incremental approach used in WRFDA. Consequently, constraints must be imposed on the lightning data to restrain large perturbations in the temperature increments from appearing in the simplified WRF adjoint model.

The other two schemes that we propose (the 1D+nD-VAR (n=3,4) approaches) assimilate lightning flash rates. Marecal and Mahfouf (2002, 2003) investigated the 1D+4D-VAR data assimilation algorithms. They used the 1D-VAR technique to adjust rainfall rates from moist physics processes (mass flux convection scheme and large scale condensation) to be closer to observed values (see Mahfouf et al. 2005). However, our assimilation procedure considers the 1D-VAR temperature retrievals to be new observations that are assimilated into the WRF model using the nD-VAR (n=3,4) WRFDA system. This approach minimizes the problem that nonlinearities in the moist convective scheme may introduce discontinuities in the cost function between the inner and outer loops of the incremental nD-VAR (n=3,4) minimization.

The present paper is structured as follows. Section 2 briefly describes the nonlinear observation operator and the lightning data used in our experiments. Section 3 introduces the WRFDA - 4D-VAR framework and describes the 3D-VAR direct assimilation approach. This includes a detailed discussion about the background error covariance matrix, observation errors, and the related covariance matrices. The final part of Section 3 provides information about the 1D+nD-VAR (n=3,4) techniques. These approaches are more consistent with the incremental WRFDA system and obviate the necessity of imposing the observation filter used in the direct assimilation approach. Sections 4 and 5 describe our use of the WRF model, the corresponding physical packages adopted in the numerical experiments, as well as the two synoptic cases. Finally, Section 6 presents results of numerical tests for both the tornadic supercell and squall line cases. Improvements in total precipitation scores, thermodynamic soundings, radar reflectivity, and convective precipitation forecasts illustrate the efficiency of the above mentioned variational assimilation techniques.

## 2. Observation operator and lightning data

The method of lightning data assimilation is highly dependent on the horizontal resolution of the meteorological model being used, in our case WRF. Barthe et al. (2010) investigated the potential for several model parameters to be used as a proxy for total lightning in cloud resolving models. Their results showed that maximum updraft velocity serves as a good proxy for flash rate in severe storms. Thus, we exploit the strong correlation between maximum vertical velocity and total flash rate given by:

$$\mathbf{f} = 5 \cdot 10^{-7} \mathbf{w}_{\mathbf{max}}^{\mathbf{k}}, \tag{1}$$

where $f$ is the total flash rate and $w_{max}$ the maximum vertical velocity in the storm. Exponent k was derived empirically to be 4.55 for continental deep convection (Price and Rind 1992). Pure parcel theory suggests that the peak vertical velocity of an updraft is related to CAPE by $w_{max} = \sqrt{2\text{CAPE}}$. However this serves only as a theoretical estimate, and Kirkpatrick et al. (2009) examined 200 convective storm simulations to relate $w_{max}$ and CAPE

$$w_{max} = 0.677 \cdot \sqrt{2\text{CAPE}} - 17.286 \tag{2}$$

Using (1-2), we link the maximum vertical velocity at each grid point to the lightning flash rate and then translate it to temperature lapse rate using CAPE. Consequently, our lightning operator assumes the form

$$\mathbf{H(X)} = 5 \cdot 10^{-7} (0.677\sqrt{2 \cdot \mathbf{CAPE(X)}} - 17.286)^{4.55}, \tag{3}$$

providing a simulated flash rate for each grid point. The input X consists of one dimensional vertical arrays of pressure, temperature, water vapor mixing ratio, and geopotential height.

Our assimilation techniques use the variational approaches described in the next section to adjust the vertical temperature profile at each grid point where innovation vectors are positive. In this way, if the model simulated lightning via CAPE is greater than observed non-zero flash rate, the scheme does not produce any direct change, keeping the same lapse rate estimated by the WRF model. Thus, we do not suppress convection where the model exhibits deep convection but observations indicate no lightning.

ENTLN observed total flash rates were used to simulate the upcoming GOES-R GLM data. The lightning data assimilated into WRF were obtained by adding the ENTLN flash rates over 20 min intervals centered on the assimilation hours and mapped on grid cells of 9 × 9 km with units of flashes $(9\ km)^{-2}\ min^{-1}$. The lightning observations for both storms being examined were generated hourly during the 2 to 6 h assimilation window between 1800 – 2000(0000) UTC, depending on the variational assimilation scheme employed for assimilation. Their relation to model calculated CAPE at locations having positive innovation vectors is shown in Fig. 1.

Equation (1) reflects the spatial 9 km resolution used to simulate the two cases and represents an adjustment to the relationship proposed by Price and Rind (1992) and Barthe et al. (2010) for cloud resolving models. Specifically, the model calculated flash rate was decreased by one order of magnitude to prevent the generation of unrealistic flash rates at large values of CAPE and to better fit the observed and calculated lightning observations (Fig. 2). For example, Fehr et al. (2004) scaled the equation by a factor of 0.26, Salzman et al. (2008) used a factor of 0.06, and Barthe and Barth 2008 used a factor of 0.19. The linear correlation coefficients between lightning observations and model simulated CAPE in Fig. 1 are 0.693 and 0.6719 for the two cases, respectively. However, if we correlate only lightning observations greater than 5 flashes $(9\ km)^{-2}\ min^{-1}$ with the corresponding CAPE generated by the model, we obtain very poor results, 0.212 for the April 2011 case and 0.124 for the June 2011 event, respectively. This indicates that the placement of deep convection by the simulation does not agree well with observations, leading to a very challenging data assimilation problem.

# 3. Methodology for lightning assimilation

## 3.1 WRFDA 4D-VAR data assimilation system

The 4D-VAR data assimilation method produces an "optimal" estimate of the true atmospheric state at the analysis time through the iterative minimization of a prescribed cost function written as

$$J(X_0) = \frac{1}{2}(X_0 - X_0^b)^T \mathbf{B}^{-1}(X_0 - X_0^b) + \frac{1}{2}\sum_{k=1}^{N}(H_k M_k(X_0) - Y_0^k)^T \mathbf{R}_k^{-1}(H_k M_k(X_0) - Y_0^k),$$

(4)

where subscript k denotes the model time step and $M_k(X_0)$ represents the state of the atmosphere at time k following k + 1 time steps of integration by the forecast model

$$M_k(X_0) = M_{0\to k}(X_0). \quad (5)$$

The 4D-VAR variational problem seeks to minimize the cost function (4) to find the analysis state $X_0$ that minimizes $J(X_0)$. This solution represents the a posteriori maximum likelihood estimate of the true atmospheric state given two sources of a priori data: the background (previous forecast) $X_0^b$ and observations $Y_0^k$ (Lorenc 1986). $H_k$ is the nonlinear observation operator that converts the initial model state into observed equivalents at time k for comparison with the corresponding observations $Y_0^k$, $\mathbf{R}_k$ is the error covariance matrix of the observations, while matrix **B** contains the background error covariances for each atmospheric variable.

The Mesoscale and Microscale Meteorology (MMM) Division of the National Center for Atmospheric Research (NCAR) currently maintains and supports the WRFDA system. The WRFDA 4D-VAR subset adopts the incremental variational formulation proposed by Courtier et al. (1994). Commonly used in operational NWP, the incremental approach is designed to find the analysis increment $\delta x = X - X_0^b$ that minimizes a cost function defined as a function of the analysis increment. Using the corresponding tangent linear observation operator and the tangent linear of the forecast model we obtain

$$J(\delta x) = \frac{1}{2}\delta x^T \mathbf{B}^{-1}\delta x + \frac{1}{2}\sum_{k=1}^{N}(d_k - \mathbf{H}_k \mathbf{M}_k \delta x)^T \mathbf{R}_k^{-1}(d_k - \mathbf{H}_k \mathbf{M}_k \delta x), \qquad (6)$$

where $d_k = Y_0^k - H_k M_k X_0^b$ are the innovation vectors at timestep k and $\mathbf{M}_k$, and $\mathbf{H}_k$ denote the tangent linear versions of the forecast model and observation operator. Two minimization algorithms are presently implemented in the WRFDA system, the conjugate gradient and the Lanczos algorithms (LaMacchia and Odlyzko 1991) which require gradient information provided by

$$\nabla_{\delta x}J = \mathbf{B}^{-1}\delta x + \sum_{k=1}^{N}[\mathbf{M}_k^T \mathbf{H}_k^T \mathbf{R}_k^{-1}(M_k(X_0) - H_k M_k(X_0^b) - \mathbf{H}_k \mathbf{M}_k \delta x], \qquad (7)$$

where $\mathbf{M}_k^T$ and $\mathbf{H}_k^T$ denote the adjoint models of the forecast model and observation operator, respectively. For NWP models, the dimension of the square symmetric and positive definite matrix $\mathbf{B}$ is $\cong 10^8 - 10^9$. Therefore, the WRFDA system does not calculate the background error covariance matrix in model space but uses a control variable v related to model space via the control variable transform U, i.e.,

$$\mathbf{B} = \mathbf{U}^T\mathbf{U} \text{ and } \delta x = \mathbf{U}v. \qquad (8)$$

Here v is the analysis increment in control variable space that includes stream function, unbalanced velocity potential, unbalanced temperature, relative humidity, and unbalanced surface pressure. To assimilate cloud and precipitation affected observations, the latest version of WRFDA not only includes the usual control variables, but also cloud water, rainwater, and ice provided by suitable parameterizations (Huang et al. 2012). Zupanski et al. (2011a) applied the maximum likelihood ensemble filter to assimilate synthetic GOES-R radiances during cloudy conditions into the WRF model. They used as control variables the potential temperature, specific humidity, and five hydrometeors (cloud water, cloud ice, rain, snow and graupel).

The observation operator described in Section 2 was implemented and integrated into the WRFDA system. The process followed the standard pattern required by WRFDA for adding a new observation operator: (1) define a new structure for manipulating and storing the lightning data; (2) define an input/output procedure adapted to the new type of observation; (3) implement quality control tests for selecting and adjusting the appropriate observations depending on the assimilation strategy; (4) implement the observation operator and develop its tangent linear and adjoint code; (5) couple the observation operator components to the minimization procedure already existing in WRFDA (conjugate gradient and Lanczos (La Macchia and Odlyzko 1991) algorithms); (6) develop verification tests for both the tangent linear and adjoint models; (7) include an observation error formulation; and (8) implement diagnostic and statistical tests. For 1D-VAR stage an additional module was created in WRFDA.

Three variational data assimilation schemes were investigated to assimilate the lightning information and evaluate the impact of the flash observations on cloud analyses, precipitation, and short term forecasts of the April 2011 and June 2011 severe storm cases.

## 3.2 Direct assimilation of lightning using WRFDA 3D-VAR

The first scheme uses the incremental three dimensional data assimilation approach to directly assimilate the lightning information. Barker et al. (2004) describe the practical implementation of a 3D-VAR system developed for the MM5 model, the precursor of the WRF model. This variational approach encounters serious difficulties when assimilating precipitation observations due to its requirements of linearity and Gaussian distributions of the observation operator and observational errors (Lopez and Bauer 2007). Highly nonlinear moist processes such as convective regimes are included in the definition of the lightning observation operator since it first converts the model variables to CAPE and then to equivalent flash observed quantities. To calculate CAPE at each observation location, vertical profiles of pressure, temperature, water vapor mixing ratio, and geopotential height are supplied. Then, the observation operator gradient information adjusts only the temperature profile in an attempt to get the model flash rates closer to the observations. The algorithm can only perform where there is at least a small amount of CAPE in the model background (otherwise the lightning sensitivities are close to zero) and only at points where the background CAPE is greater than or equal to 325.973 J kg$^{-1}$. This value was calculated via observation operator formula (3) since the bracketed term must be positive to generate real positive flash rates. The performance of the variational system depends on the background error statistics. These statistics contain important information about how an observation spreads in model space and how the final analysis is physically balanced. The background error covariance matrix **B** describes the probability density function of the forecast errors. We estimated **B** using ensemble statistics. Using this method, vertical and horizontal error covariances are represented by empirical orthogonal functions and a recursive filter.

The lightning observations were assumed to be uncorrelated. Consequently, the observation error covariance matrix is diagonal. For simplicity our assimilation tests were performed using an identity matrix for the observation error covariance matrices **R**.

As mentioned by Zupanski et al. (1993), Zou (1997), and Vukicevic and Bao (1998), observation operators that contain highly nonlinear processes as in our case require special treatments such as thresholds, and switches must be applied to avoid discontinuities. Otherwise the analysis may not be optimal, or even worse, the cost functional minimization can fail. A measure of how well the incremental variational assimilation system will perform is how linearly the observation operator behaves.

The linearity test used in this study is derived from the alpha test described by Navon et al. (1992). It is based on comparing the output of the tangent-linear model with those from finite difference calculations using the forward model. The tangent linear model was obtained using the TAPENADE Automatic Differentiation Engine (Hascoët and Pascual 2004). A function of α

$$F(\alpha) = \frac{H(x+\alpha\delta x) - H(x)}{\alpha \mathbf{H}\delta x} \tag{9}$$

can be used as a measure of linearity in which $\delta x$ represents the initial perturbation of the control vector and α is a scaling factor. For small values of α that are not too close to the double precision machine epsilon, one expects F(α) to be close to 1 for linear behavior. However, δx variations must be carefully analyzed since δx theoretically could be so small that even nonlinear models may exhibit a nearly linear performance. An appropriate choice for perturbations would be the analysis minus first guess departures from the 3D-VAR analysis. These differences give an indication of the maximum perturbations that can be encountered during the minimization procedure. Therefore, we performed the linearity test for $\delta x = x_a - x_b$ which is scaled with α ranging from $10^{-9}$ to 1. Since the innovation vectors are decreased by changing only the temperature profiles, we tested how linearly the lightning operator behaves with respect to temperature perturbations.

Fig. 3 depicts plots of $\log_{10} |1 - F(\alpha)|$ from five differently located vertical profiles of the WRFDA 3D-VAR model analysis at 1800 UTC 27 April 2011. The legend box displays the magnitude of the innovation vectors responsible for temperature and CAPE perturbations that occurred in the analysis. If $\log_{10} |1 - F(\alpha)| \leq -1$, a reasonably linear regime is achieved. For example, for $\log_{10} |1 - F(\alpha)| = -1$, F(α) is between 0.9 and 1.1, i.e., a 10% disagreement between the finite difference and the tangent linear models.

For very small values of α (Fig. 3), the results deteriorate due to numerical round off. For large values of α, the results are clearly affected since larger innovation vectors produce larger temperature perturbations and inflict stronger nonlinearities. For α close to 1 and for innovation vectors greater than 2 (see green line in Fig. 3), the lightning observation operator ceases to behave linearly ($\log_{10}|1 - F(\alpha)| \geq -1$), leading to inconsistency between the inner and outer loop calculations of the cost functional in the incremental variational WRFDA system. Consequently, sub-optimal analyses are obtained, or even worse, the minimization algorithm may fail. Several strategies can be adopted to avoid these difficulties, including restrictions on temperature perturbations inside the minimization algorithms or imposing a threshold for the innovation vectors participating in the data assimilation process. The second approach is supported by the results presented in Fig. 3, where magnitudes of temperature perturbations are strongly sensitive to the input innovation vectors. Fortunately, the above results show the analysis minus first guess perturbations of temperature and CAPE and not the intermediary variations obtained during the minimization iterations which usually are much smaller. Consequently a more relaxed threshold can be used, and in one of our 3D-VAR experiments, only the lightning observations with innovation vectors less than or equal to ten flashes $(9 \text{ km})^{-2} \text{ min}^{-1}$ were assimilated. The results are discussed in Section 6.

We next propose more sophisticated strategies consisting of 1D+nD-VAR (n=3,4) data assimilation schemes that we implemented into the WRFDA system.

## 3.3 Data assimilation of lightning using WRFDA 1D+3D-VAR and 1D + 4D-VAR

The 1D+4D-VAR approach was first used by Gérard and Saunders (1999) to assimilate Special Sensor Microwave Imager (SSM/I) clear sky radiances in the ECMWF model. Their methodology was adapted and implemented by Marecal and Mahfouf (2002) in the ECMWF 4D-VAR system (Rabier et al. (2000) to assimilate surface rain-rate retrievals from SSM/I and the TRMM Microwave Imager. The two step technique became operational at ECMWF during June 2005 (Bauer et al. 2006a, b) and is currently utilized for assimilating SSM/I microwave brightness temperatures in cloudy and rainy regions. Lopez and Bauer (2007) studied the potential impact of assimilating NCEP stage IV analyses of hourly accumulated surface precipitation over the U.S mainland using the 1D+4D-VAR method.

The Japan Meteorological Agency (JMA) has implemented a 1D+4D-VAR technique for radar reflectivity observations to improve analyses of water vapor and precipitation using the non-hydrostatic mesoscale JNoVA model (Honda et al. 2005). One dimensional (1D) pseudo observations of relative humidity (RH) are retrieved from radar reflectivity, with the retrieved RH data assimilated as conventional data in JNoVA. Caumont et al. (2010) reported that the resulting precipitation forecasts were improved by combining the 1D retrieval and a three-dimensional variational (3D-VAR) data assimilation method (1D+3D-VAR).

To avoid nonlinearity effects introduced by our lightning operator in the direct assimilation scheme, and to extract the full benefit from all available flash observations, we used the 1D-VAR+nD-VAR (n=3,4) technique similar to that outlined by Marecal and Mahfouf (2002). For the remainder of the present paper, when we refer to nD-VAR techniques, we mean 3D&4D-VAR variational data assimilation approaches. In the first step of the procedure, the raw lightning measurements are used in the 1D-VAR technique to produce increments of temperature that are added to the model background to generate column temperature retrievals in accordance with the flash observations. Then, these temperature pseudo observations are assimilated as conventional observations into the WRFDA nD-VAR systems.

The 1D-VAR method searches for an optimal estimate of the temperature profile at the analysis time through the iterative solution of the following prescribed cost function,

$$J(X_0) = \tfrac{1}{2}(X_0 - X_0^b)^T \mathbf{B}^{-1}(X_0 - X_0^b) + \tfrac{1}{2}(\tfrac{H(X_0)-y_0}{\sigma_0})^2, \tag{10}$$

where H is the lightning observation operator (3) introduced in the previous section and $\sigma_0$ is the observation variance. The 1D-VAR retrieval method is equivalent to a best estimate approach based on Bayes theorem. The 1D-VAR problem basically is solved in a similar way to 3D-VAR, except that only the vertical dimension is considered.

The National Meteorological Center (NMC) method proposed by Parrish and Derber (1992) was applied to construct the background error vertical covariance matrix B for the temperature profiles. We used 12 h and 24 h forecasts valid at the same time from a one month dataset generated by the WRF model.

Figure 4 presents correlation coefficients between average temperatures at different altitudes at locations of lightning observations for the two storm cases being studied. Strong correlation is evident between vertical levels 10 and 25 (1,500 - 6,000 m) for the supercell storm case on 27 April 2011. An even stronger correlation is seen between levels 25 and 35 (6,000 - 10,000 m) for the 15 June 2011 squall line case.

For the 1D-VAR variational retrieval, we employed the CONMIN minimization algorithm proposed by Shanno and Phua (1980) to represent the quasi-Newton limited memory conjugate gradient. A survey of four different conjugate gradient methods for large scale minimization problems in meteorology can be found in Navon and Legler (1987). The control vector $X_0$ in our application contains vertical profiles of temperature on 60 model levels. The minimization requires the gradient of $J(X_0)$ defined in (11):

$$\nabla J(X_0) = \mathbf{B}^{-1}(X_0 - X_0^b) + \mathbf{H}^T \left(\tfrac{H(X_0)-y_0}{\sigma_0}\right), \tag{11}$$

where the adjoint of the observation operator $\mathbf{H}^T$ again was calculated using TAPENADE (Hascoët and Pascual (2004)). The 1D+nD-VAR methods have several advantages. They introduce additional quality control tests, better handle the less linear inversion problem, present 'smooth' pseudo observations to nD-VAR, filter the other nD-VAR control variables, and use B twice. However, these approaches are more computationally expensive than traditional nD-VAR data assimilation schemes.

## 4. NWP model

We selected the non-hydrostatic Weather Research and Forecasting (WRF) model version V3.3 with the advanced research dynamical core (ARW) for the lightning data assimilation experiments (Wicker and Skamarock 2002, Skamarock et al. 2005, 2008). WRF-ARW was configured (Fig. 5) with an outer domain with 27 km horizontal grid spacing and a 9 km horizontal grid spacing covering a two way nested inner domain of approximately 1413 km $\times$ 1170 km for both storm events. 60 vertical levels were selected to cover the troposphere and lower part of the stratosphere between the surface to approximately 20 km. The grid size of the 9km model domain is $157 \times 130 \times 60$. For initial and boundary conditions the NCEP Global Forecasts System (GFS) 1 degree resolution final analyses were used. In terms of model physics, the Kain-Fritsch cumulus parameterization scheme (Kain and Fritsch (1993) was utilized, and the Yonsei planetary boundary layer scheme (Hong and Dudhia 2003) was chosen. For the radiative transfer of long- and short-wave radiation we selected, respectively, the rapid radiative transfer model (RRTM) (Mlawer et al. 1997) and the Dudhia (1989) scheme. The RRTM spectral scheme interacts with resolved clouds and accounts for multiple bands, trace gases, and microphysics species, while the Dudhia scheme uses simple downward integration and includes efficient water vapor, cloud albedo, and clear sky absorption and scattering. For microphysical processes a single moment, 6 class, cloud microphysics scheme (Hong and Lim 2006) was used. It includes graupel and has demonstrated good performance in high resolution simulations of the amount of precipitation and temporal evolution of observed cases of heavy rainfall.

## 5. Synoptic cases

We ran simulations for two severe thunderstorm cases over the United States. The cases included the extensive outbreak of deadly supercell tornadoes on 27 April 2011 and a squall line causing numerous severe wind reports across the Southeast on 15 June 2011. These two days were among the most electrically active for the eastern U.S. during 2011.

Fig. 5 depicts the geographical domain of study which covers the Southeast U.S. All of the simulations included 6 h of model spin up between 1200 UTC and 1800 UTC, after which lightning assimilation began with an assimilation window varying between 2 to 6 h. The simulations then were run an additional 3 - 7 h without assimilation, ending at 0300 UTC of the next day.

## 6. Results

We examined the three proposed variational data assimilation schemes described above: 3D-VAR, 1D+3D-VAR, and 1D+4D-VAR for each of the two storm cases. Two experiments were performed with each assimilation scheme. The control variables in WRFDA V3.3 include stream function, unbalanced potential velocity, unbalanced temperature, unbalanced surface pressure, and pseudo relative humidity. We first selected only unbalanced temperature as the control variable (configuration I – C1) to measure the impact of lightning data assimilation. The second experiment was a sensitivity test in which additional control variables were added, stream function, unbalanced velocity potential, unbalanced surface pressure, and pseudo relative humidity (configuration II – C2), to decrease the initial innovation vectors. For the 3D-VAR scheme we performed an additional experiment that limited the number of assimilated lightning observations according to the linearity test results discussed in Section 3.2 (configuration I + QC).

In the case of the 3D-VAR and 1D+3D-VAR schemes, a cycling procedure was adopted to assimilate the lightning observations between 1800 UTC and 0000 UTC. Seven hourly 3D-VAR assimilations were used during the assimilation window. The first guesses were obtained by integrating the previous 3D-VAR analysis 1 h in time using the WRF model. For the 1D+4D-VAR scheme we used a 2 h assimilation window between 1800-2000 UTC. We also attempted to assimilate the 1D-VAR temperature profiles derived from flash rates over a 6 h time interval using 4D-VAR, but the results were not successful since both of the minimization algorithms implemented in WRFDA failed to converge. Specifically, the gradient test exhibited poor results. The radiation and cumulus parameterizations are computationally expensive to call every time step; so the recommendation is to run these routines depending on the spatial resolution of the model. At 9 km resolution the parameterizations must be called each nine time steps, and their simplified adjoint models lead to large errors in the gradient calculation.

The number of lightning observations utilized by the 3D-VAR scheme was greater than those used by the 1D+nD-VAR schemes. This is explained by the number of quality control (QC) tests imposed on the lightning observations during the 1D-VAR stage. The vertical background error covariance matrix acted as a virtual QC test that conditioned the 1D-VAR minimization. Figure 6 shows an example of the average reduction of both the cost function and its gradient during the 35 successful 1D-VAR retrievals using data from 1800 UTC 27 April 2011. A logarithmic scale is used. The cost function gradients are reduced by two orders of magnitude for the observation components. Figure 6 also reveals that most retrievals require 80 iterations to achieve a one order of magnitude decay in the cost functional.

Since our knowledge of lightning observation errors still is limited, we chose to eliminate from the assimilation process those observations that only led to a small decrease in values of the innovation vectors. The numbers of successful 1D-VAR retrievals and rejected observations is shown in Fig. 7 for both storm cases. The number of observations rejected due to 1D-VAR convergence failure is greater during the 27 April event. One possible explanation is the differing flash rates. The assimilated average flash rates for the 27 April and 15 June storms were 1.317 and 2.378 flashes $(9 \text{ km})^{-2} \text{ min}^{-1}$, respectively. Thus, the $H(x_a) - H(x_b) < 0.2$ condition allowed the nD-VAR procedures to assimilate profiles of temperature obtained from areas where deep convection occurred according to the intensity of lightning. The average flash rates used to generate the successful 1D-VAR retrievals were 5.343 and 7.065 flashes $(9 \text{ km})^{-2} \text{ min}^{-1}$, respectively, for the 27 April and 15 June storms.

To increase the performance of the 1D+nD-VAR WRFDA lightning assimilation schemes, we applied an additional quality control test for proxy temperatures derived via the 1D-VAR assimilation system. Specifically, the retrievals were checked for vertical consistency (super adiabatic lapse rate) and were adjusted to dry adiabatic in unstable layers using OBSPROC (the WRFDA Running Observation Preprocessor).

The number of lightning observations directly assimilated by the 3D-VAR scheme can be seen in Fig. 7. However, one should note that the accepted observations include the flash rates rejected during the 1D-VAR phase (shown as red and blue). Thus, the restriction strategies imposed in 1D-VAR are successful since more accurate precipitation results are obtained from the 1D+3D-VAR scheme than the direct 3D-VAR approach during the assimilation window.

Average temperature increments from the 1D-VAR CONMIN algorithm (Fig. 8) reveal warming near the surface that decreases the stability (increases the CAPE) and thereby strengthens convection at locations where lightning is observed. On 27 April the atmosphere between the surface and ~ 600 m is warmed to increase the convection intensity and better match the observed flash rate, while on 15 June the warming affects a deeper layer up to an average of 1200 m. Thus, our strategy is producing the desired results that are consistent with parcel theory. Specifically, if a parcel near the surface becomes warmer than its environment, it becomes buoyant and is more likely to reach its level of free convection (LFC), form a cloud, and possibly produce lightning.

Figures 9-10 depict the efficiency of the 1D-VAR assimilation scheme, showing values of the innovation vectors before (left panels) and after (right panels) ingesting the lightning observations. The 1D-VAR minimization algorithm decreases the maximum value of the innovation vectors from 11.53 to 2.91 flashes $(9\ km)^{-2}\ min^{-1}$ and from 20.15 to 3.16 flashes $(9\ km)^{-2}\ min^{-1}$ for the 27 April and 15 June storms, respectively. Thus, at the observation locations, we found new vertical temperature patterns that increase CAPE and generate flash densities closer to the observed ones following the lightning observation operator formula.

The following parts of this section describe results from assimilating lightning using the 3D-VAR scheme as well as assimilating the 1D-VAR temperature pseudo observations for the 1D+nD-VAR schemes using only temperature as the control variable. Figures 11 and 13 (left panels) show locations of lightning observations along with the initial differences between model calculated flash rates and the observed rates used by the 3D-VAR direct approach at 1800 UTC on both 27 April and 15 June. During the April event, the most intense lightning activity was located in northeastern Alabama. Other locations of lightning include Arkansas, Mississippi, Missouri, North Carolina and Tennessee. The 3D-VAR scheme modifies the background temperatures to increase the amount of CAPE in these areas in accordance with the flash rates (right panels). The maximum increment of ~ 1200 $J\ kg^{-1}$ occurs near the border between northeast Alabama and Tennessee. The 3D-VAR minimization algorithm fails to acknowledge the lightning in southeastern Missouri and central Tennessee due to the small CAPE in the background state (see discussion in Section 3.2).

The greatest CAPE increment occurs at 34.66° N, 86.14° W in northeast Alabama, with Fig. 12 showing the direct effects of lightning assimilation at that location. The greatest changes occur in the 850-1000 hPa and 150-200 hPa layers. The atmosphere is warmed near the surface, with increases as great as 3-4° C. Consequently the level of free convection (LFC) is decreased from ~ 850 hPa to 1000 hPa, while the equilibrium level (EL) is increased from 200 hPa to 170 hPa. The vertical profile of dew point temperature is not modified.

During the June squall line case (Fig. 13, left), areas of lightning are located mostly in Tennessee, with sparse flashes in South Carolina and Florida. The direct 3D-VAR approach increases the background temperatures, thereby producing greater CAPE increments with values exceeding 1400 J kg$^{-1}$ over Tennessee and Florida (Fig. 13, right) where the greatest innovation vectors occur. Locations 35.93° N, 86.92° W and 28.74° N, 82.46° W are where the innovation vectors exhibit the greatest values of 20.15 and 15.13 flashes (9 km)$^{-2}$ min$^{-1}$. The greatest change in CAPE occurs at 28.85° N, 82.83° W, near the location with the second greatest innovation vector. The profiles of temperature and dew point temperature before and after 3D-VAR assimilation of lightning are illustrated in Fig. 14. Altitudes of the LFC and EL are modified in agreement with the increased CAPE.

The 1D-VAR temperature profile retrievals are used as observations in the 1D+nD-VAR assimilation schemes. Locations of these pseudo observations are shown in Figs. 9 and 10 at 1800 UTC for both storm cases. In the case of the 27 April tornado outbreak, an isolated observation in North Carolina (Fig. 9) allows us to examine the impact of the background error covariance matrix on the 1D+3D-VAR analysis. The intensity of lightning is 3.3 flashes (9 km)$^{-2}$ min$^{-1}$, and maximum generated CAPE increments are 400 J kg$^{-1}$ in the vicinity of the lightning observation (Fig. 15, left). The background error covariance matrix spreads the lightning information across North Carolina, Virginia, and South Carolina, increasing the amount of CAPE in the 1D+3D-VAR analysis proportional to the distance from the lightning observation. Figure 15 shows the utility of 1D+3D-VAR since the greatest change in CAPE is in northeastern Alabama where the most intense flashes are located. Western Mississippi is a secondary region where the presence of lightning has increased the CAPE.

For the 15 June case at 1800 UTC, the greatest increments of CAPE generated by the 1D+3D-VAR lightning data assimilation method are across Tennessee, South Carolina and Florida (Fig. 16, left), locations of the 1D-VAR retrievals (Fig. 10).

A comparison between CAPE increments obtained by the 3D-VAR (Fig. 11, 13 − right panels) and 1D+3D-VAR assimilation schemes (Fig. 15, 16 − left panels) reveals that greater values are initiated by the 3D-VAR approach. The impact of the assimilation with respect to total precipitation is discussed in the last part of this section.

The 1D+4D-VAR increments of CAPE are depicted in Figs. 15 and 16 (center and right panels) for both cases. The center panels in Figs. 15 and 16 were obtained by employing only temperature as the control variable (configuration I − C1), while the right panels show results when four additional control variables were used: stream function, velocity potential, surface pressure, and relative humidity (configuration II − C2). Figure 17 reveals that greater increments of CAPE are generated using C2. The 27 April case is more sensitive to the control variable configurations in terms of CAPE increment magnitudes. The changes in CAPE patterns at analysis time 1800 UTC are explained by the locations of lightning observations, their flash rates, and the nature of the 4D-VAR algorithm. The algorithm uses the WRF model as a strong constraint to seek the initial condition that best fits the 1D-VAR pseudo observations within the entire 2 h assimilation interval between 1800 and 2000 UTC. The following experiments were performed using the C1 configuration.

We next assess how the WRF model propagates the analyses generated by the lightning assimilation methods. We compared the analyses utilizing the lightning observations with a control simulation having no assimilation. NCEP GFS 1 degree resolution final analyses were utilized to derive the initial and boundary conditions at 1200 UTC on both days, followed by a 6 h spin up period. The minimization packages contained in the 3D-VAR and 1D+nD-VAR schemes utilized the control state as the first guess at the start of the assimilation window.

Figure 18 compares simulated radar reflectivities generated by the control, 3D-VAR, and 1D+3D-VAR simulations 10 min after the analysis time (1800 UTC) for the 27 April case. The CAPE perturbations across Georgia (Fig. 11 − right panel and Fig. 15 left panel) slightly reduce rain water mixing ratio between 930 and 730 hPa immediately after assimilation (1801 UTC, not shown). Consequently, an area of reflectivity of 5-10 dBZ, denoting light rain, forms over Georgia in both assimilation runs between 1801 and 1820 UTC (Fig. 18). These weak reflectivities do not suggest lightning, but they are a result of the assimilation schemes. The area of the 3D-VAR generated precipitation (center panel) is slightly larger than the one produced by the 1D+3D-VAR simulation (right panel). This is explained by the number of lightning observations assimilated by each scheme. Specifically, 425 observations were assimilated using 3D-VAR approach, while 1D+3D-VAR assimilated 35 observations.

Unlike 3D-VAR analyses, the 4D-VAR analyses initialize the reflectivity fields to zero. Thus, we compare the 1D+4D-VAR and control reflectivities at the end of 1D+4D-VAR assimilation window (2000 UTC) to allow a spin up period of 2 h. Figure 19 shows the potential of the 1D+4D-VAR scheme for lightning assimilation. Reflectivity structures of moderate precipitation (15-35 dBZ, right panel) possibly capable of producing scattered lightning, are seen across eastern Tennessee where two small areas of lightning are observed (Fig. 19 − left panel). These reflectivities are considerably greater and more contiguous than those from the control (center panel). Flash observations and reflectivity areas also are well correlated over northeastern Arkansas and Mississippi. The theoretical aspects of the 4D-VAR confer an important advantage over the 3D-VAR counterpart that is seen in Fig. 19. The effect of lightning in the 4D-VAR simulation occurs at the same time with the observed lightning, not after a spin up period as in the 3D-VAR case. The 4D-VAR scheme suffers from this impediment only at the beginning of the assimilation window

Simulated radar reflectivity plots from the control, 3D-VAR, and 1D+nD-VAR simulations for the 15 June case at 2010 UTC are presented in Fig. 20. Of the various assimilation procedures, the 1D+4D-VAR scheme produces the greatest changes in reflectivity compared to the control. Reflectivity areas associated with heavy rain are observed in the 1D+4D-VAR simulation across eastern Tennessee where lightning is observed. Areas of enhanced reflectivity also are observed in the 1D+4D-VAR run over southeastern Kentucky, West Virginia and Virginia. However, they are poorly correlated with the observed lightning at 2000 UTC. Areas of precipitation in western North Carolina are reduced by both the 3D-VAR and 1D+3D-VAR schemes, more pronounced in the former case and are consistent with the 1D+4D-VAR reflectivity pattern.

We next examine how efficiently 1D+4D-VAR initiates and places precipitation and how accurately the amounts compare with output from the 3D-VAR and 1D+3D-VAR simulations. Specifically, we consider 1 h precipitation totals between 1900 and 2000 UTC on both study days (Figs. 21 and 22). The simulated totals are verified against hourly precipitation from the stage IV dataset. Stage IV data are determined from radar and rain gauge observations and benefit from human quality control at the NWS River Forecast Centers (Lin and Mitchell 2005) . The approximately 4 km stage IV data were interpolated to the 9 km grid spacing of our simulations using the average cell grid method.

We first discuss the 27 April case (Fig. 21). The 1D+4D-VAR simulated rainfall totals are improved over the control in Arkansas, Missouri, southern Illinois, Tennessee and Virginia.

However, over Alabama and Mississippi 1D+4D-VAR produces rain areas that are too large compared with stage IV. Although the control simulation produces areas of moderate to heavy rain totals over southeastern Missouri, southern Illinois, and southwestern Virginia, 1D+4D-VAR suggests smaller totals that are more consistent with the stage IV information. Also, the large area of moderate precipitation over central Mississippi in the control simulation is much weaker in the 1D+4D-VAR simulation, better matching stage IV. Another effect of 1D+4D-VAR assimilation is seen over North Carolina and South Carolina where areas of control precipitation are reduced. The 3D-VAR and 1D+3D-VAR simulations show similar precipitation patterns and increase the size of the precipitation area near the border between northeastern Alabama and Tennessee where lightning occurs at 1800 and 1900 UTC.

Hourly precipitation ending at 2000 UTC 15 June is illustrated in Fig. 22. Once again, 1D+4D-VAR reproduces the stage IV precipitation fields much better than the other simulations. The rain areas in Kentucky, southern Illinois, Mississippi, and western Tennessee are reduced in size in the 1D+4D-VAR simulation and better reproduce the stage IV patterns. The area of precipitation from 1D+4D-VAR is more contiguous than the patchy precipitation from the other simulations, better agreeing with stage IV.

We next quantify differences between the various versions of simulated precipitation and stage IV observations by calculating root mean square errors during the assimilation windows (1800 – 0000 UTC for the cycled 3D-VAR and 1D+3D-VAR assimilation schemes; 1800 – 2000 UTC for the 1D+4D-VAR approach) and the 3 to 7 h forecast period between 2000 (0000) - 0300 UTC when no assimilation was performed. Seven simulations are compared to the control simulation. The simulations (except for the control) were initialized with the 3D-VAR and 1D+nD-VAR analyses using the C1 and C2 settings, while the C1+QC configuration was used only for the direct 3D-VAR assimilation scheme.

The performances of our schemes relative to stage IV are illustrated in Fig. 23. The 1D+4D-VAR simulations exhibit the smallest errors during the assimilation window. The results are impressive since the RMSE scores are improved (compared to the control) an average of 25% (1D+4D-VAR –C2) and 14% (1D+4D-VAR –C1) on April 27. On June 15 the errors are decreased (compared to the control) an average of 27.5% (1D+4D-VAR –C2) and 12% (1D+4D-VAR –C1). Best results are obtained using the C1 setting when only temperature was used as a control variable. During the 1D+4D-VAR forecast window, i.e. 2000 - 0300 UTC, the RMS errors of 1D+4D-VAR deteriorate quickly on April 27, but on 15 June the errors continue to maintain the same level of improvement as during the assimilation window. This finding is explained by the relative accuracy of the control simulations used as the first guess. These statistical results quantify the subjective findings in Figs. 21 – 22 that the 27 April case was better forecast than 15 June.

RMSE values of the cycled 3D-VAR and 1D+3D-VAR schemes (Fig. 23) document noticeable improvements in precipitation skill compared to stage IV between 2000 and 0000 UTC on both days. Compared to the 1D+4D-VAR algorithm, the 3D-VAR approaches produce smaller changes in precipitation areas and intensities. The 3D-VAR and 1D+3D-VAR simulations need a 2 h time window to start benefiting from the lightning information; this is not needed in the 1D+4D-VAR approach. RMSEs of 1D+3D-VAR are smaller than those from the direct 3D-VAR scheme, indicating the affinity of WRFDA incremental approaches to linear observation operators. Overall the 1D+3D-VAR schemes increase precipitation accuracy by 5.33% and 9% (compared to the control) during the forecast period between 2000 and 0000 UTC on both days. RMSEs of 1D+3D-VAR are not sensitive to the choice of the control variable configuration (C1 and C2) since similar results were obtained. The 3D-VAR precipitation skill also is improved, with the RMSE reduced by the following factors on the April and June cases, respectively: 2.96% and 4.143% (C1), 2.53% and 2.49% (C2), and 2.61% and 5.31% (C1+QC). During the C1+QC experiments we imposed a threshold to limit the number of assimilated lightning observations based on the magnitude of the corresponding innovation vectors. Specifically, based on the linearity test described in Section 3.2, we eliminated lightning observations having innovation vectors greater than 10 flashes $(9\text{ km})^{-2}\text{ min}^{-1}$. The restriction assured the validity of the observation operator tangent linear assumption. The above RMSE results show that this strategy was successful only on 15 June when larger innovation vectors occurred. Magnitudes of the innovation vectors were shown in the left panels of Figs. 11 and 13.

## Conclusions

Three variational schemes, 3D-VAR and 1D+nD-VAR (n=3,4), have been developed to assimilate lightning data into the WRF numerical forecast model. Our lightning assimilation techniques used observed flash rates to warm the atmosphere near the surface, increase the CAPE, and thereby strengthen the simulated convection at locations of lightning. The observation operator used in our experiments is highly nonlinear, and is a modified version of the expression found in Barthe et al. (2010) adapted for 9 km resolution. Simulations with a horizontal grid spacing of 9 km were performed across the central and eastern United States on two of the most electrically active days during 2011.

Results showed that hourly precipitation patterns, its statistics, and radar reflectivity were improved by assimilating the lightning observations. Values of RMSE revealed that each of the proposed lightning data assimilation methods improved simulated precipitation relative to stage IV observations during the assimilation windows. The 1D+4D-VAR approach performed best, improving the precipitation areas and totals by 25% and 27.5% compared to the control run on the two days that were studied. RMSE scores of the 1D+4D-VAR simulations also were the best during a subsequent 7 h forecast period on 15 June. However, on 27 April the 1D+4D-VAR forecasts outside the assimilation window were not improved due to the relatively high accuracy of the control simulation. RMSEs of the 3D-VAR and 1D+3D-VAR forecasts did not decrease following the assimilation period.

The results indicated that lightning observations provide varying amounts of information depending on the choice of the control variables in the variational minimization algorithms. We tested two control variable configurations. The precipitation RMS errors showed that the best configuration includes only temperature as control variable during the assimilation window. Adding more control variables such as stream function, velocity potential, surface pressure, and relative humidity didn't improve the precipitation errors.

Additional storm cases must be investigated at 9 km resolution to confirm the usefulness of our lightning observation operator since its formulation was chosen independently of the storm cases that were studied. And, additional tests must be performed to decide which control variable configuration is more reliable in the context of lightning data assimilation. Certain aspects like improving observation error correlations based on instrument measuring errors and the background error covariance matrix approximations may improve the lightning assimilation.

The 1D+3D-VAR method proved to be superior to the 3D-VAR scheme, since the former used a smaller number of lightning observations and produced more accurate precipitation results. The WRFDA incremental approach requires that nonlinear observation operators behave linearly which is not valid for all observations in the assimilation process. This causes errors to accumulate in the gradient calculations. Another factor that must be considered is the accuracy of the observation error correlations. The quality control tests incorporated during the 1D-VAR stage acted as a thinning algorithm to decrease the negative effects of using erroneous correlations. Also, the 1D+nD-VAR approaches benefit twice from the background error covariance matrix.

The number of observations assimilated by the proposed methods can be increased by including observations that have negative innovation vectors and optimizing the 1D-VAR minimization algorithm to decrease the algorithm's failure rate. The 1D-VAR algorithm failed to converge mostly because of the lack of CAPE in some regions of the first guess fields. An additional step, such as a nudging scheme, that would artificially increase background CAPE to an amount that allows the observation operator to sustain the lightning data assimilation would alleviate the convergence failure rate of 1D-VAR minimization algorithm.

Results of the 1D+4D-VAR lightning assimilation and short term forecasts indicate improvements in precipitation scores and show promise for operational implementation. This will help make the upcoming GOES-R lightning data a valuable source of information over the United States and regions of poor or no radar coverage such as portions of the Atlantic and Pacific Oceans, Caribbean, Central America, and South America.

## Acknowledgments

This research was sponsored by NOOA award N A10N ES4400008 under the auspices of the Joint Center for Satellite Assimilation.

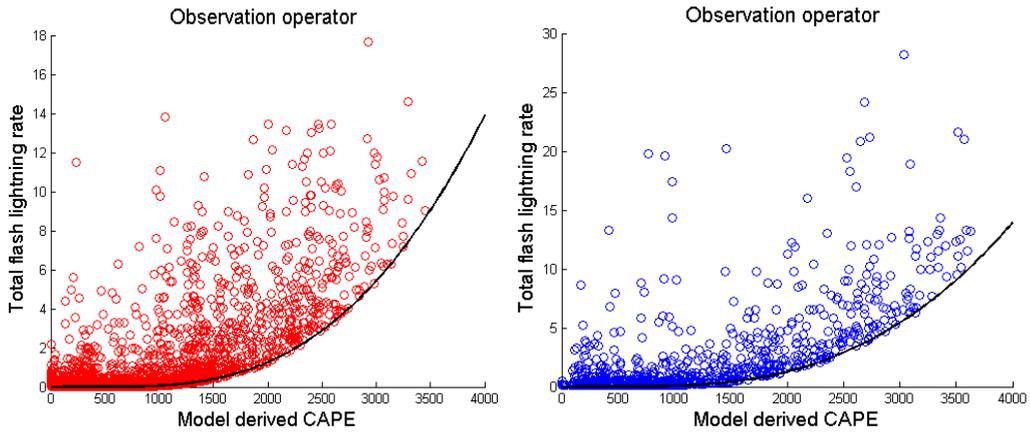

Figure 1.The distribution of lightning observations vs. model simulated CAPE for the 27 April 2011 supercell (left) and the 15 June 2011 squall line (right) cases. The black curves depict the observation operator in (3).

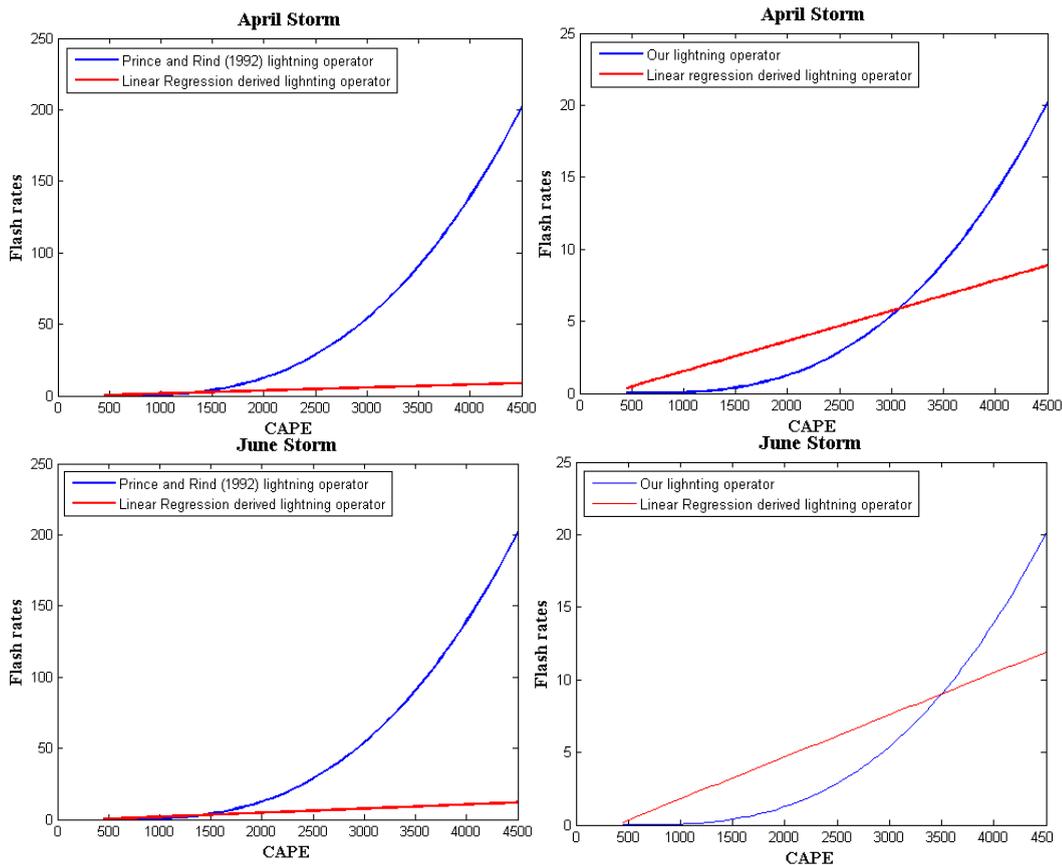

Figure 2. Comparison between the Prince and Rind (1992) lightning observation operator, linear regression empirical derived lightning observation operators, and our lightning observation operator in (3).

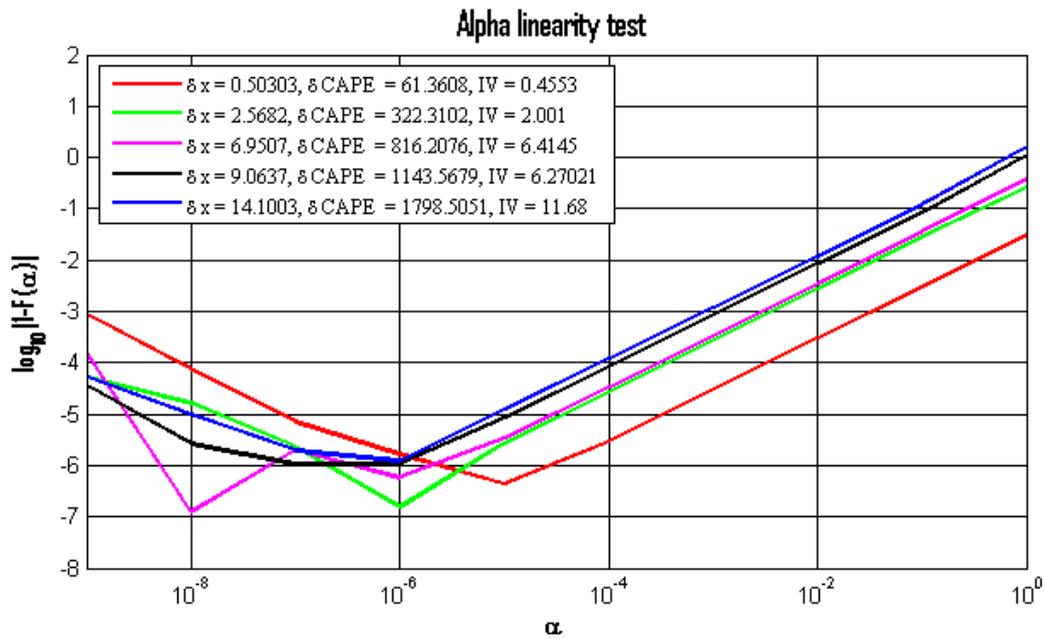

Figure 3. Linearity test for the lightning observation operator for temperature perturbations varying from 0.5 to 14.1 K, and the corresponding change in CAPE from 61.3 to 1798.5 J kg$^{-1}$. Temperature perturbations occur at each vertical level, and δx represents a general measure of it being calculated using the Euclidian norm.

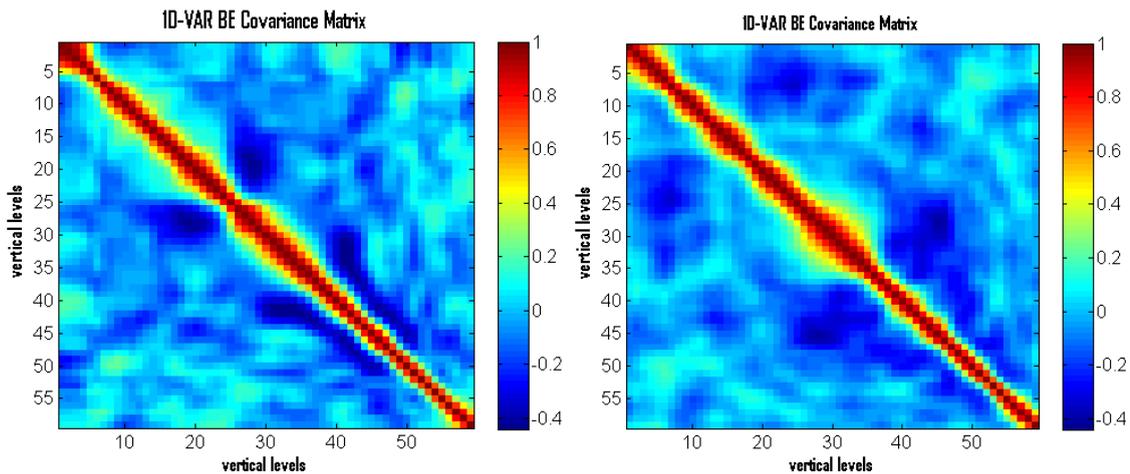

Figure 4. Average vertical correlation coefficients for temperatures at locations of lightning for the 27 April supercell (left) and 15 June squall line (right) cases.

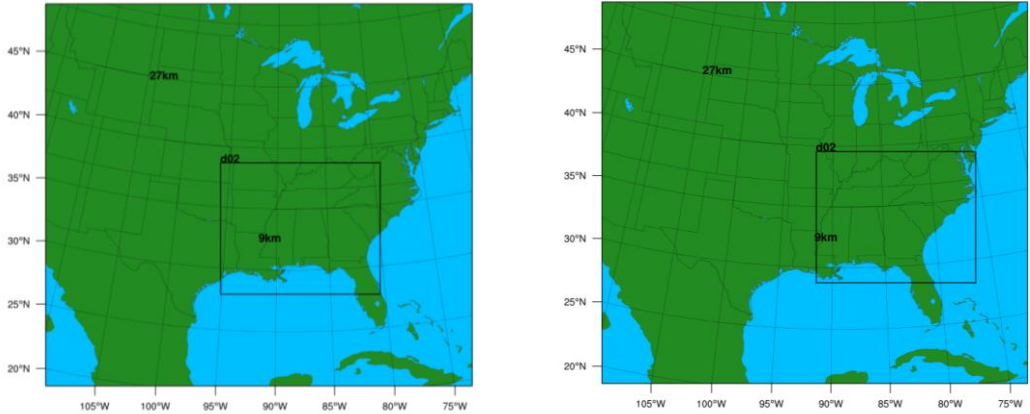

Figure 5. Computational domains for the 27 April supercell (left) and 15 June 15 squall line (right) cases.

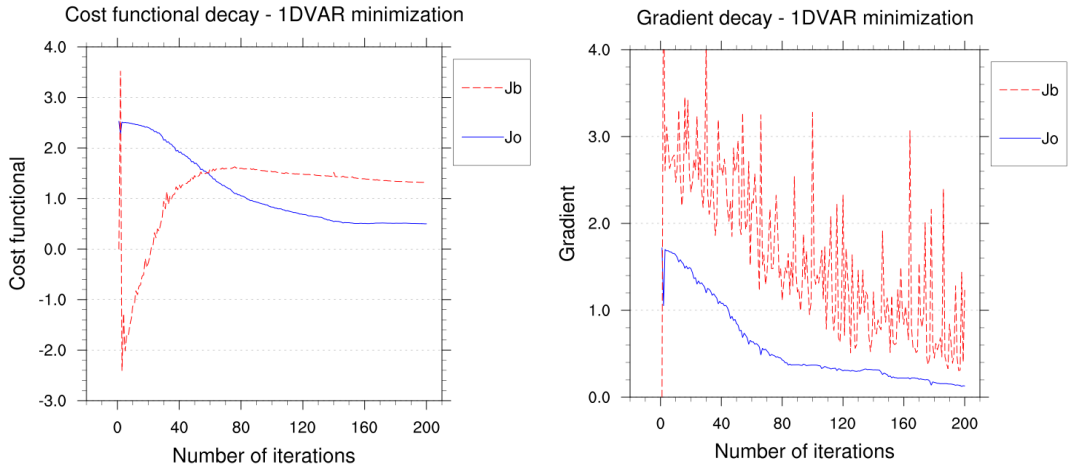

Figure 6. Mean of the 1D-VAR cost function $J$ and gradient $\nabla J$ for the observation (blue) and background (red) terms using data from 1800 UTC 27 April.

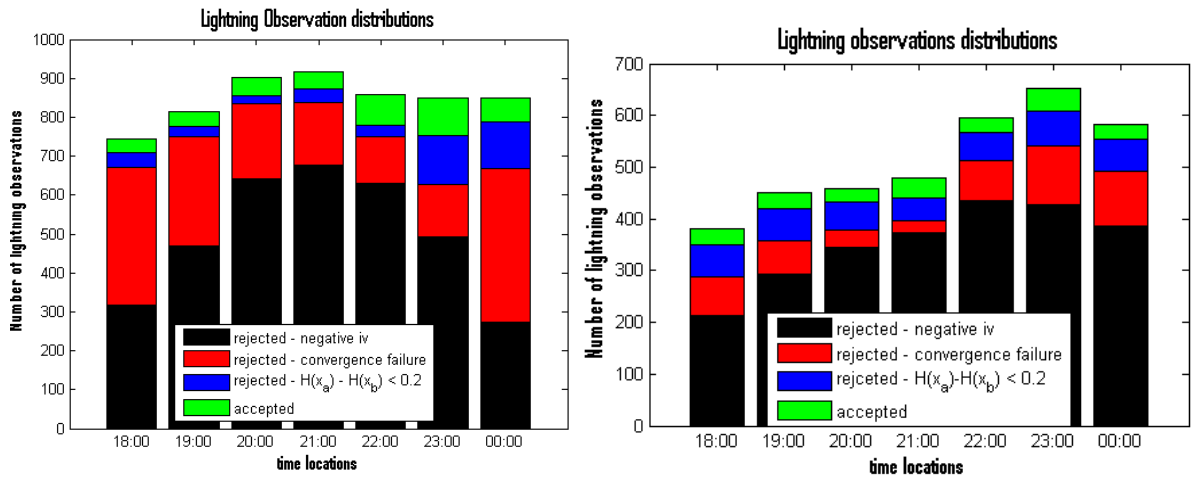

Figure 7. Histograms of rejected observations and successful 1D-VAR retrievals for the 27 April (left) and 15 June (right) cases.

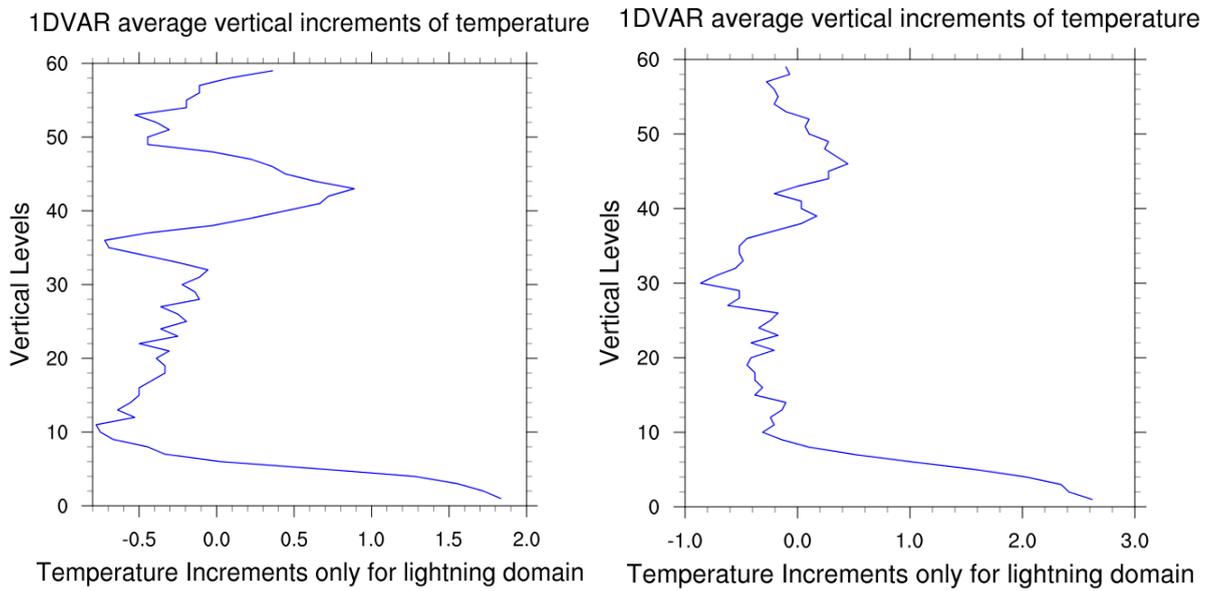

Figure 8. Average vertical increments of temperature (K) for the successful 1DVAR retrievals at 1800 UTC on 27 April (left) and 15 June (right).

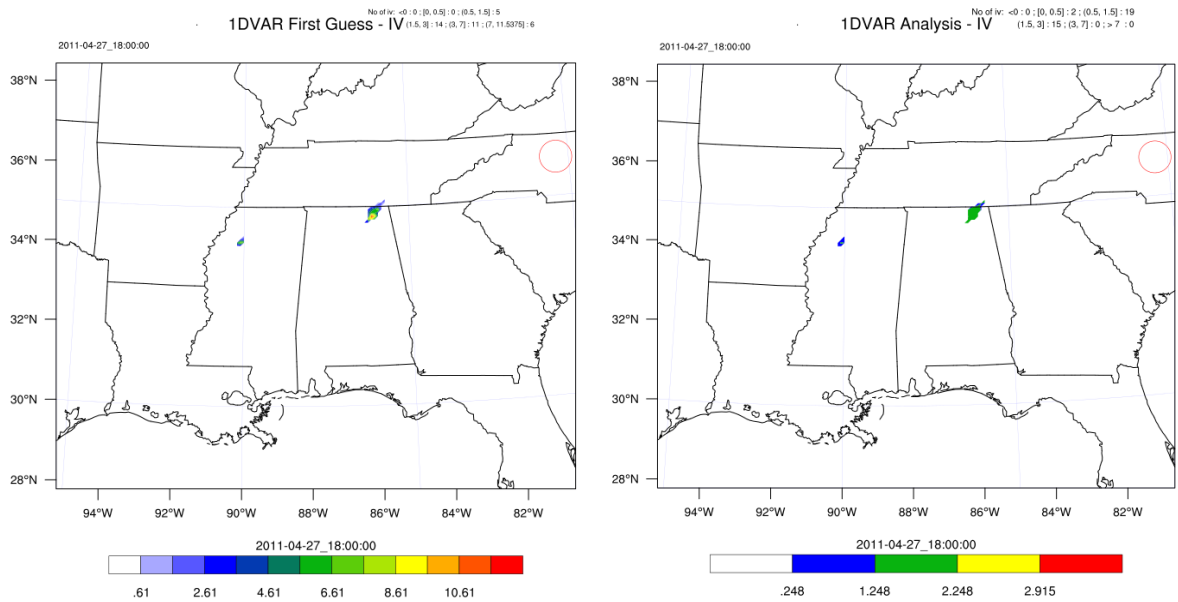

Figure 9. Innovation vectors (flashes (9 km)$^{-2}$ min$^{-1}$) before (left) and after (right) 1DVAR assimilation at 1800 UTC 27 April. Scattered (difficult to see) lightning observations are located in the center of the red circle having a radius of 0.5 $^{o}$ longitude.

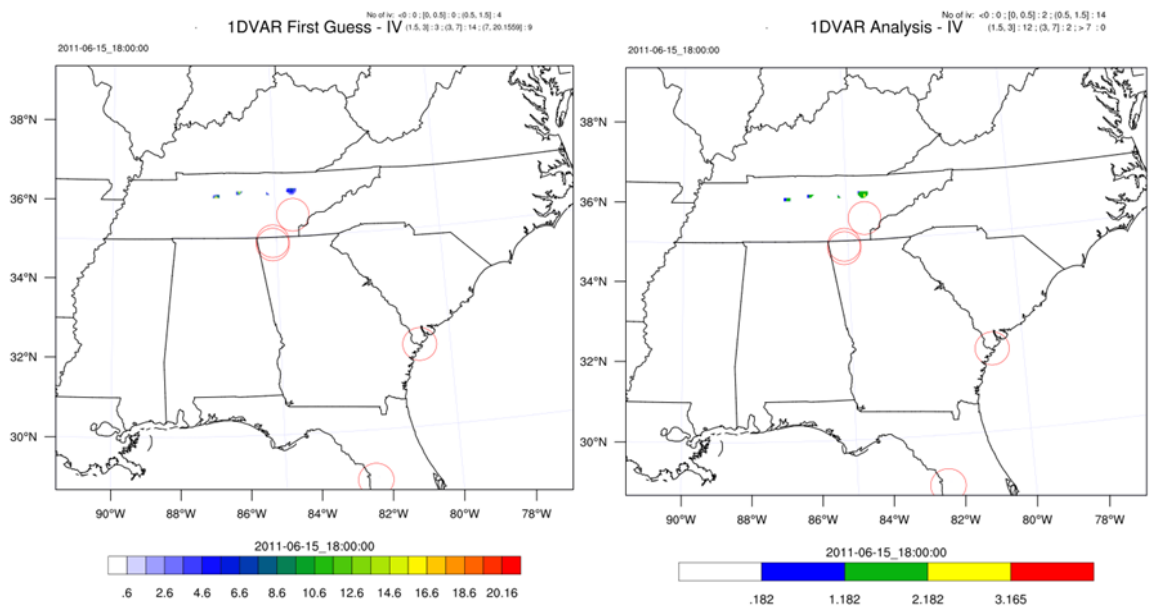

Figure 10. Innovation vectors (flashes (9 km)$^{-2}$ min$^{-1}$) before (left) and after (right) 1DVAR assimilation at 1800 UTC 15 June. Scattered (difficult to see) lightning observations are located in the center of the red circles having a radius of 0.5 $^{o}$ longitude.

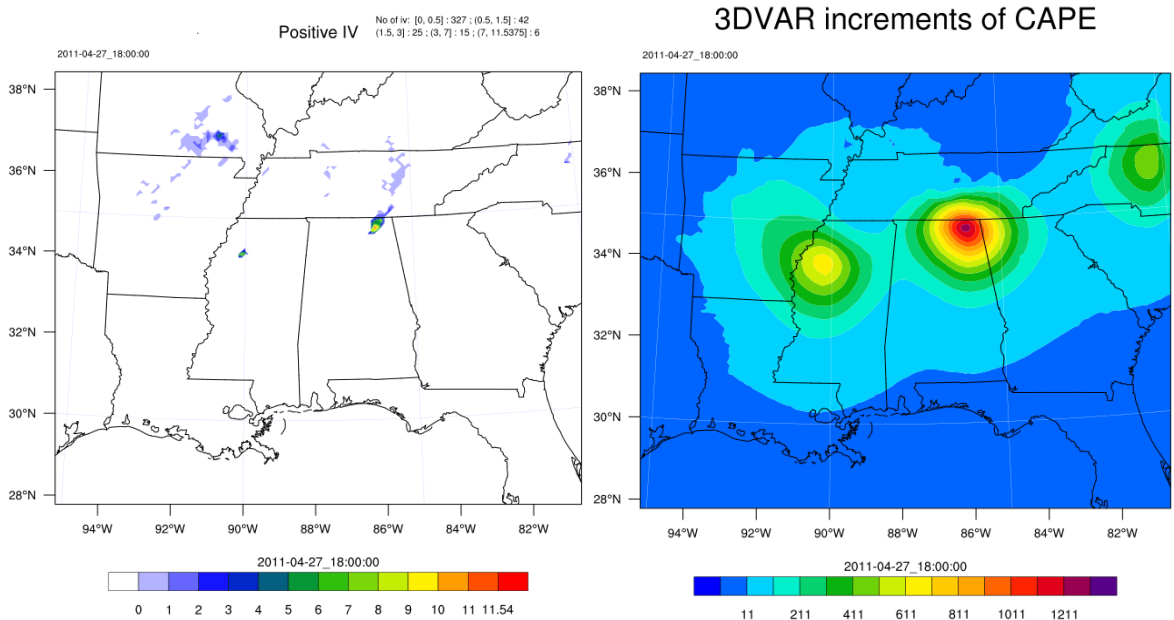

Figure 11. Innovation vectors (flashes (9 km)$^{-2}$ min$^{-1}$) before (left) and corresponding increments of CAPE (right; J kg$^{-1}$) following 3DVAR lightning assimilation at 1800 UTC 27 April 27.

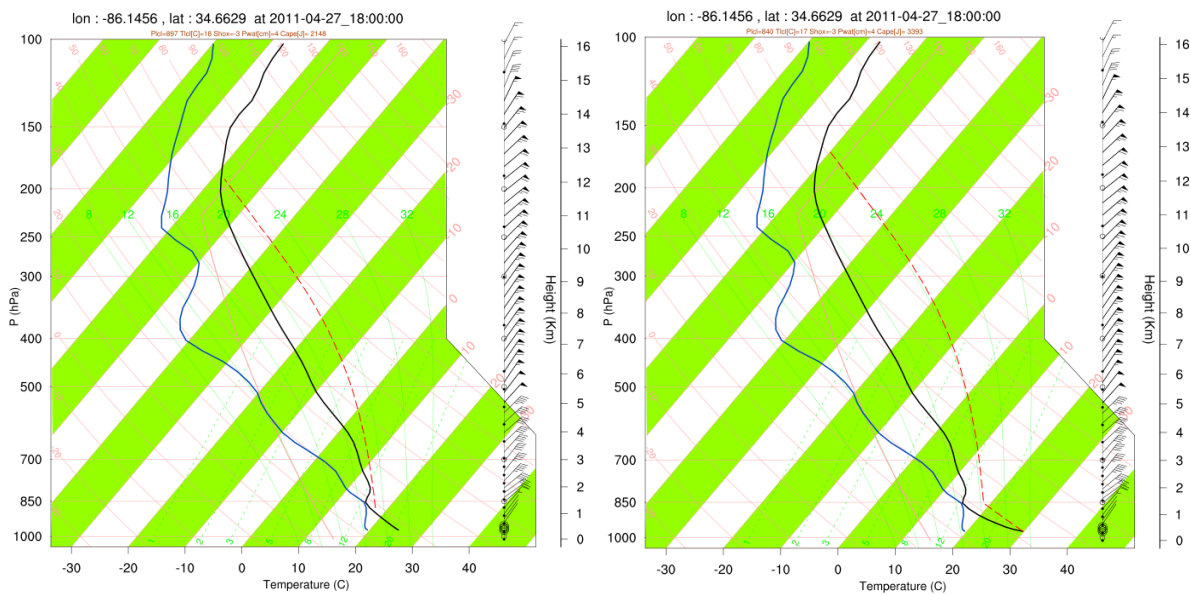

Figure 12: Skew-T diagrams (left, no lightning; right, after 3DVAR assimilation of lightning) at 1800 UTC 27 April at the location of greatest change in CAPE observed in northeast Tennessee. Air temperature (°C, black line), dew point temperature (°C, blue line), and horizontal wind (kt, barbs along right axis) are shown.

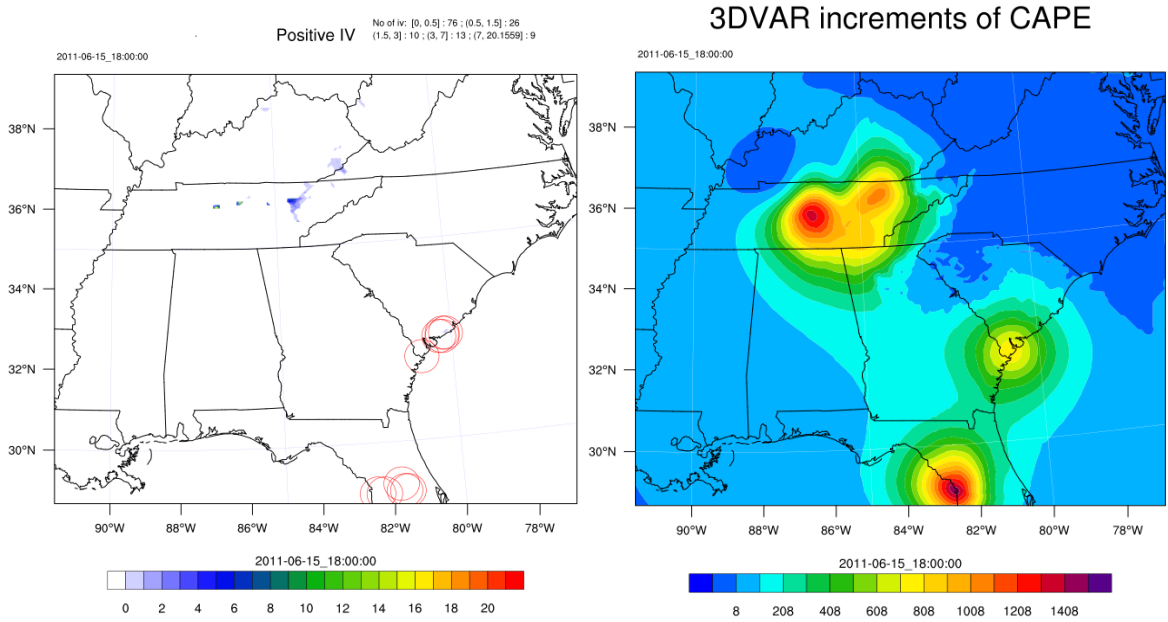

Figure 13. Innovation vectors (flashes (9 km)$^{-2}$ min$^{-1}$) before (left) and corresponding increments of CAPE (right, J kg$^{-1}$) following 3DVAR lightning assimilation at 1800 UTC 15 June.

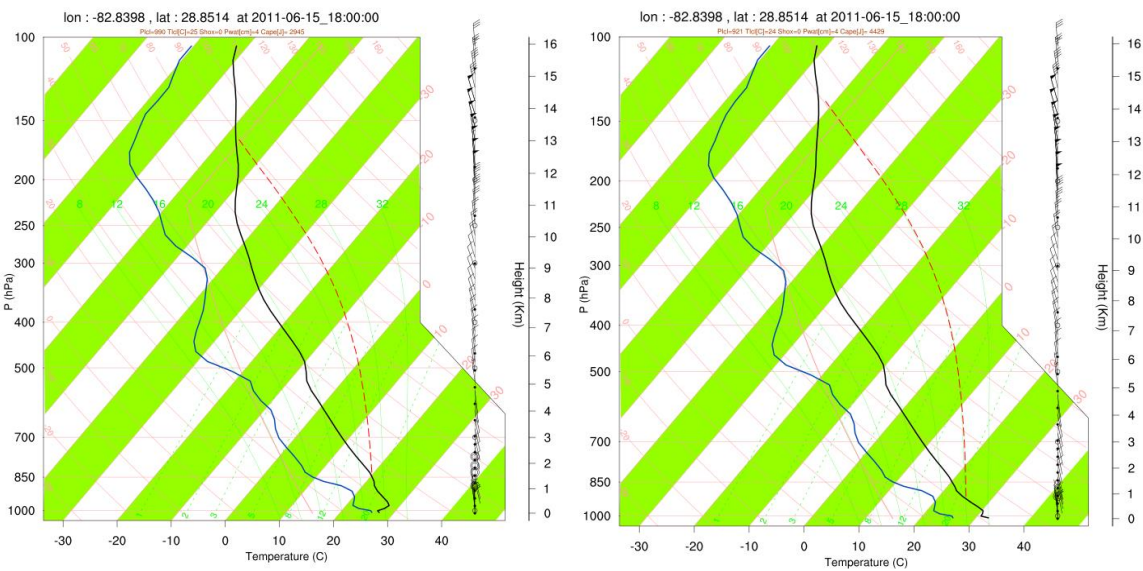

Figure 14. Skew-T diagrams (left, no lightning, right, after 3DVAR assimilation of lightning) at 1800 UTC 15 June at the location of greatest change in CAPE located in central Florida with air temperature (°C, black line), dew point temperature (°C, blue line), and horizontal wind (kt, barbs along right axis).

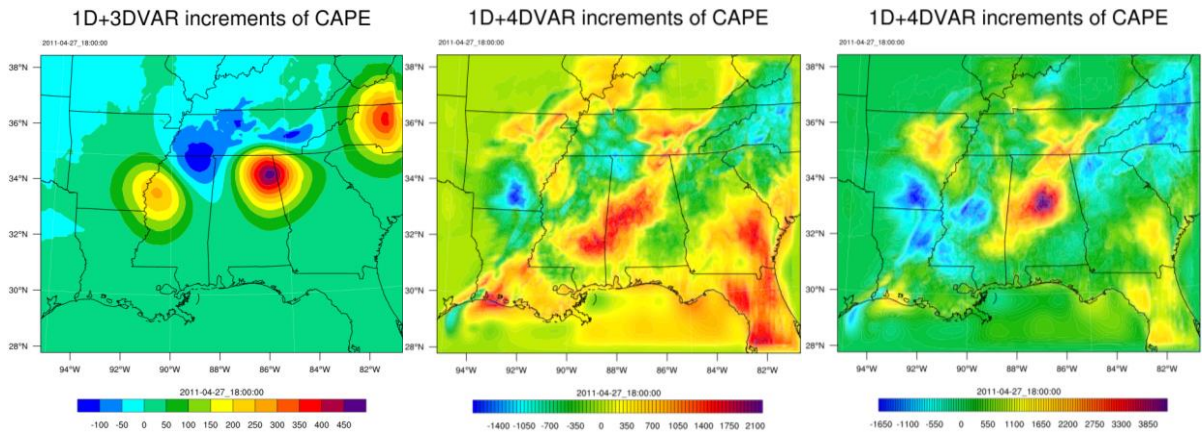

Figure 15. Increments of CAPE (J kg$^{-1}$) following 1D+3DVAR (left), 1D+4DVAR (center) with temperature as control variable, 1D+4DVAR (right) with stream function, velocity potential, temperature, surface pressure and relative humidity as control variables following the lightning assimilation at 1800 UTC 27 April.

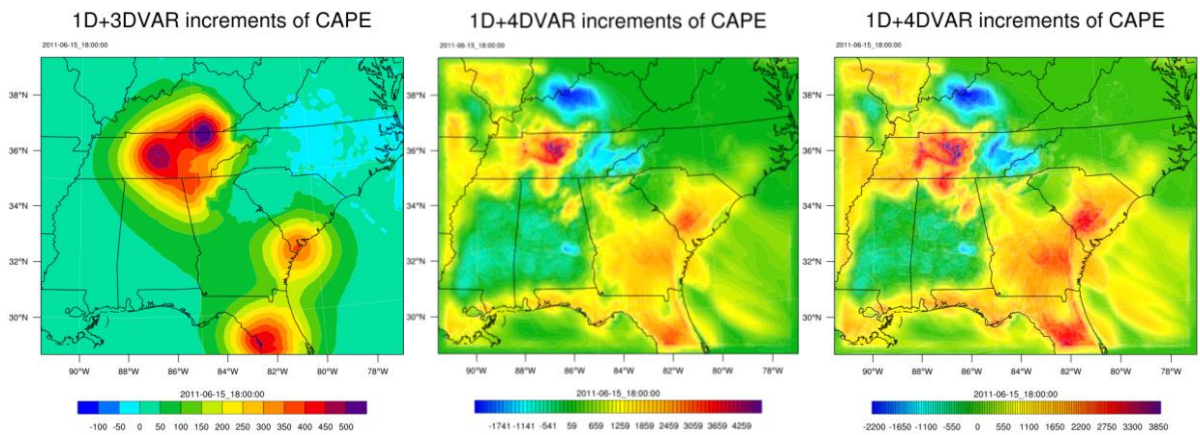

Figure 16. Increments of CAPE (J kg$^{-1}$) following 1D+3DVAR (left), 1D+4DVAR (center) with temperature as control variable, 1D+4DVAR (right) with stream function, velocity potential, temperature, surface pressure and relative humidity as control variables following the lightning assimilation at 1800 UTC 15 June.

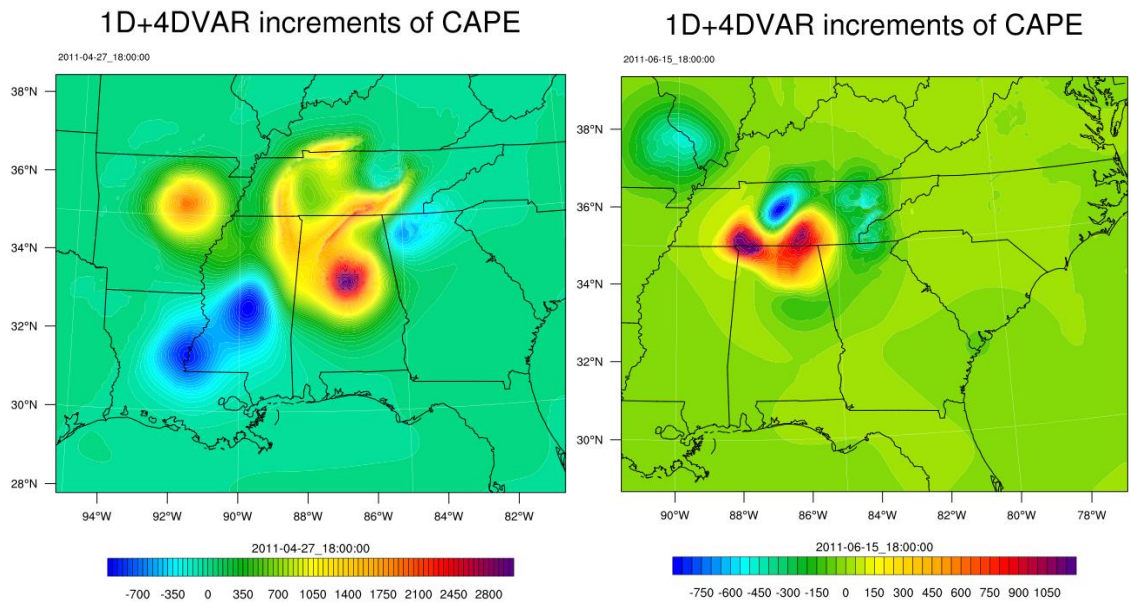

Figure 17. Differences in CAPE patterns (J kg$^{-1}$) obtained with the 1D+4DVAR strategies using five control variables (stream function, velocity potential, temperature, surface pressure and relative humidity) and one control variable (temperature) following the lightning assimilation at 1800 UTC 27 April (left) and 15 June (right).

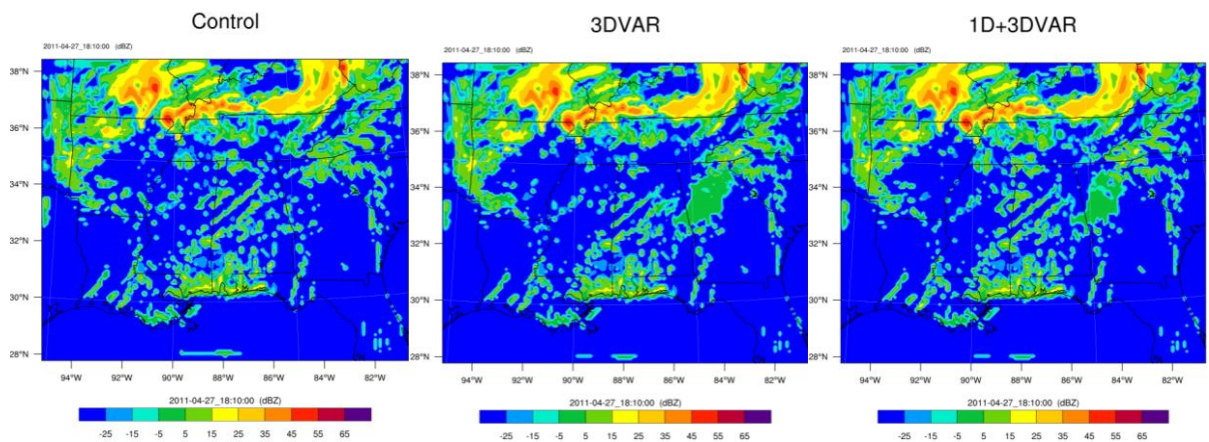

Figure 18. Simulated radar reflectivity (dBZ) at 1810 UTC 27 April from the 3D-VAR and 1D+3D-VAR procedures.

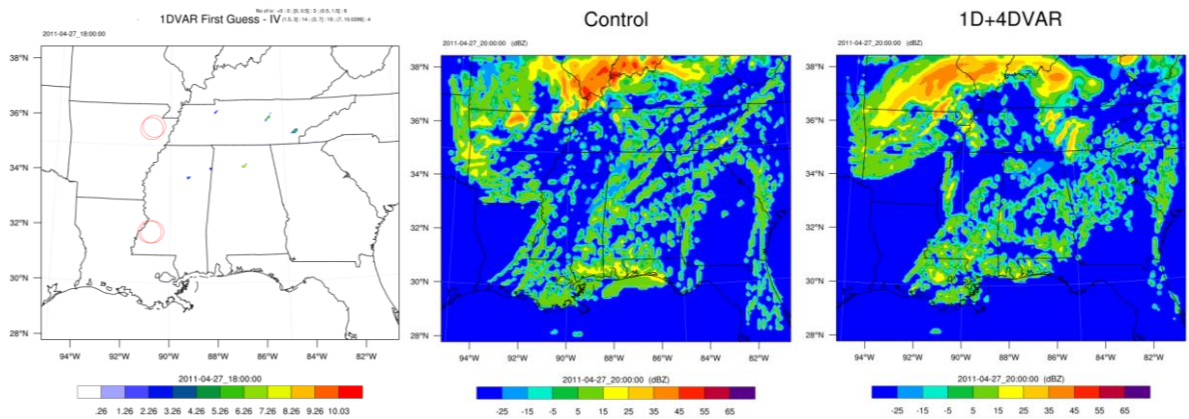

Figure 19. Simulated radar reflectivity (dBZ) at 2000 UTC 27 April from the 1D+4DVAR and control simulations.

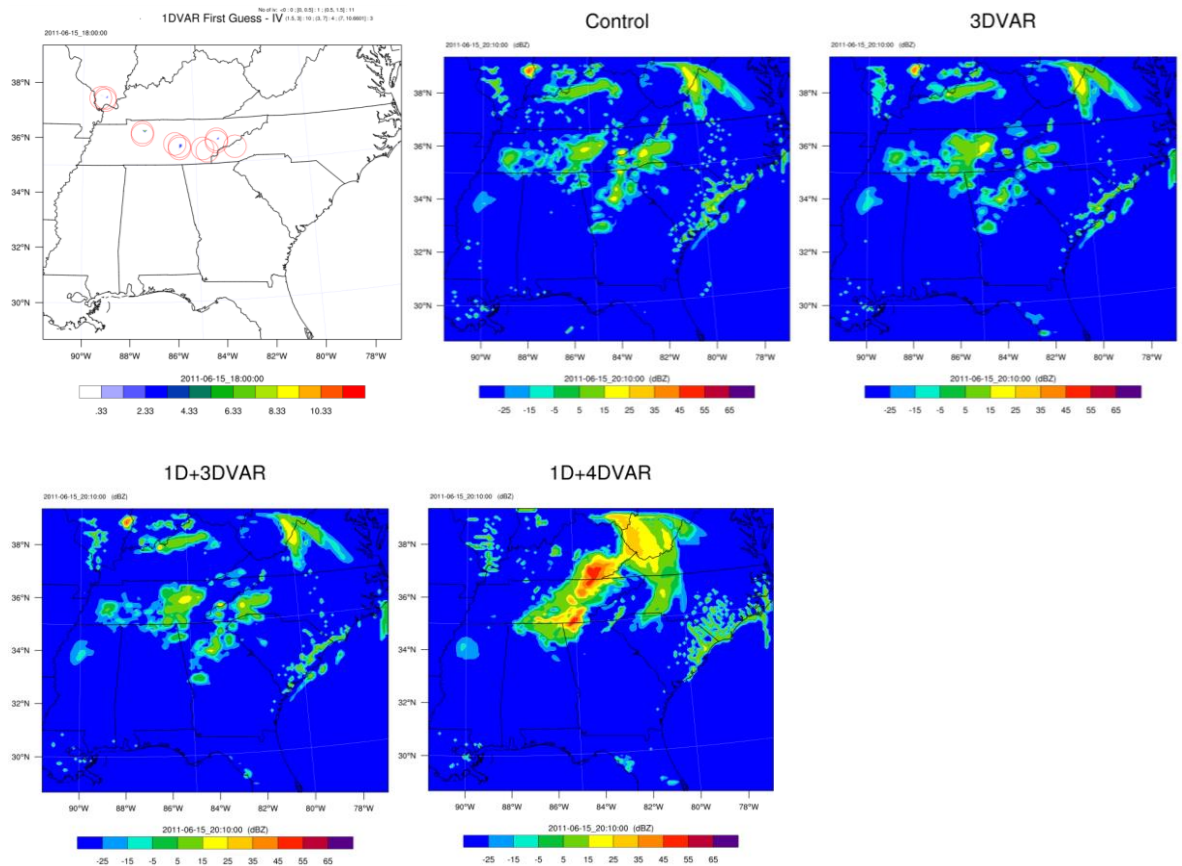

Figure 20. Simulated radar reflectivity (dBZ) at 2010 UTC 15 June from all the assimilation procedures.

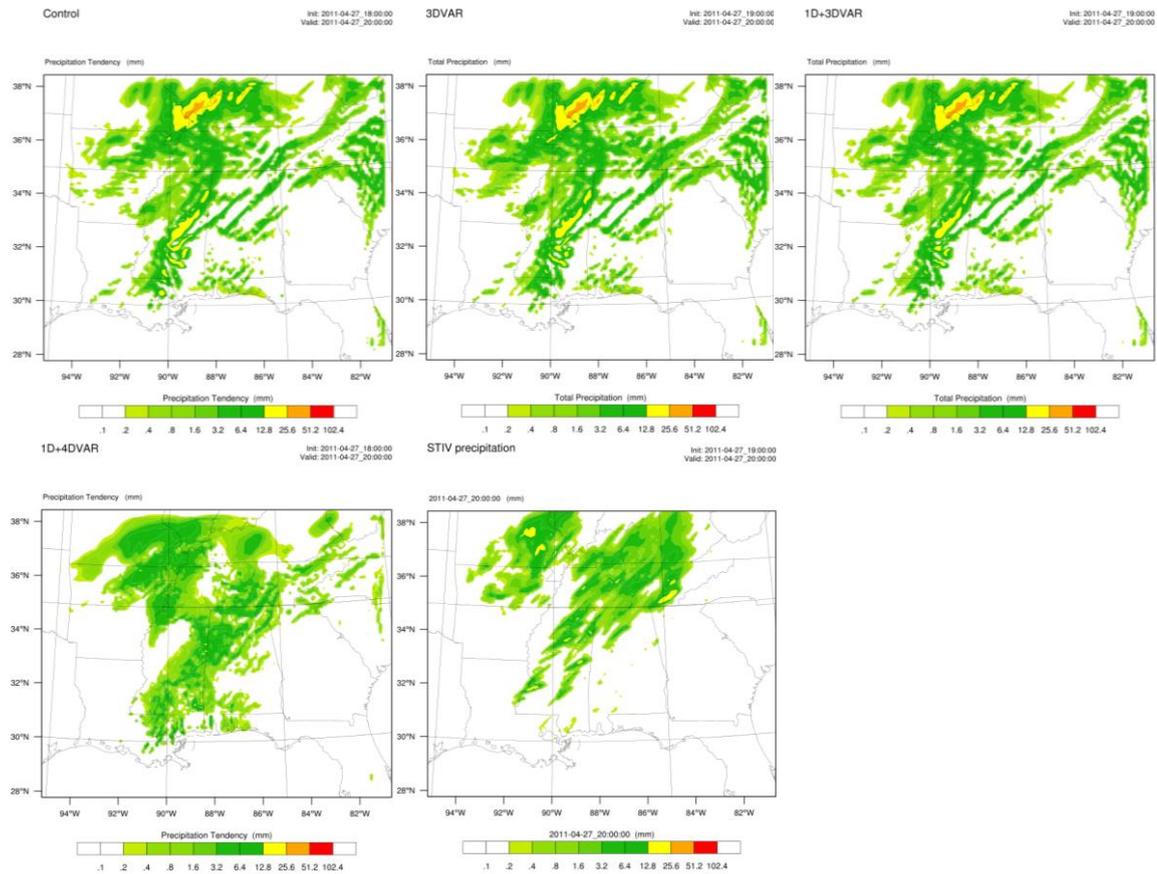

Figure 21. 1 h precipitation (mm) ending at 2000 UTC 27 April from the control run, various assimilation procedures, and stage IV observed precipitation.

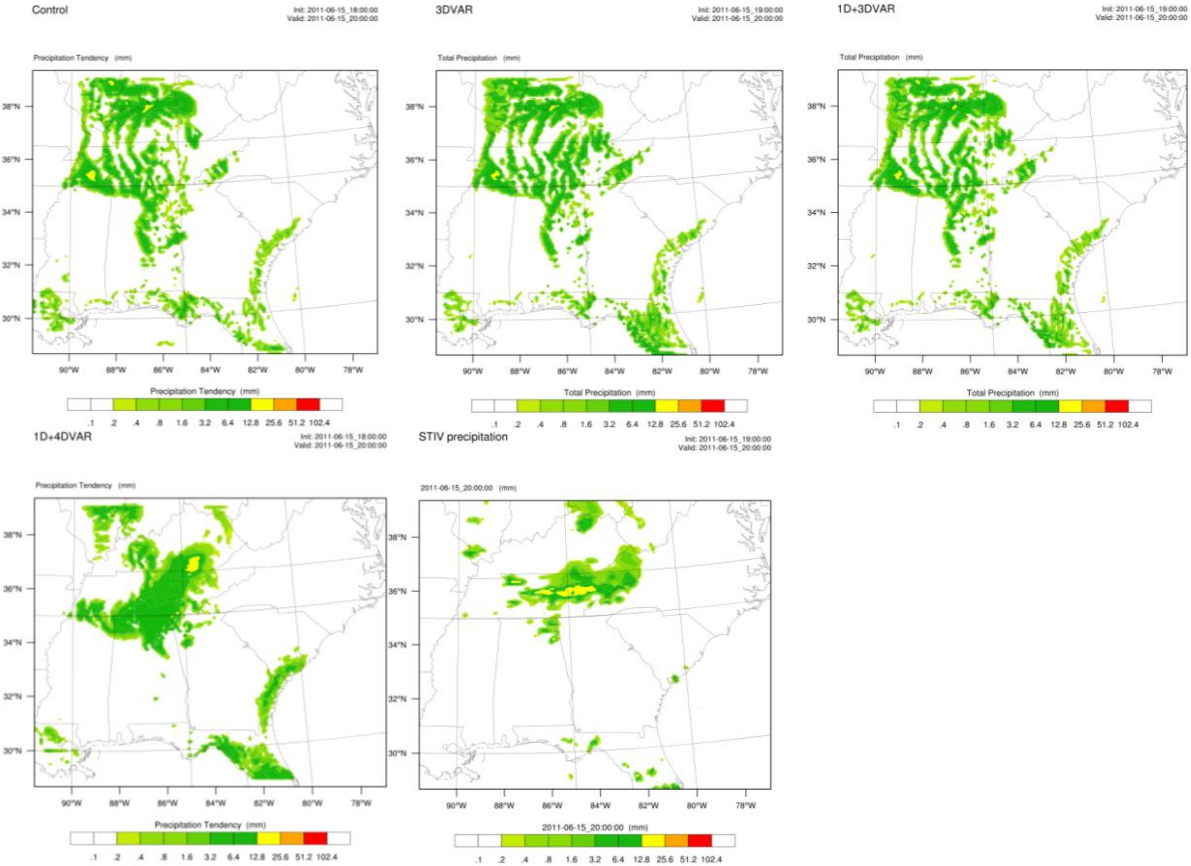

Figure 22. 1 h precipitation (mm) ending at 2000 UTC 15 June from the control run, various assimilation procedures, and stage IV observed precipitation.

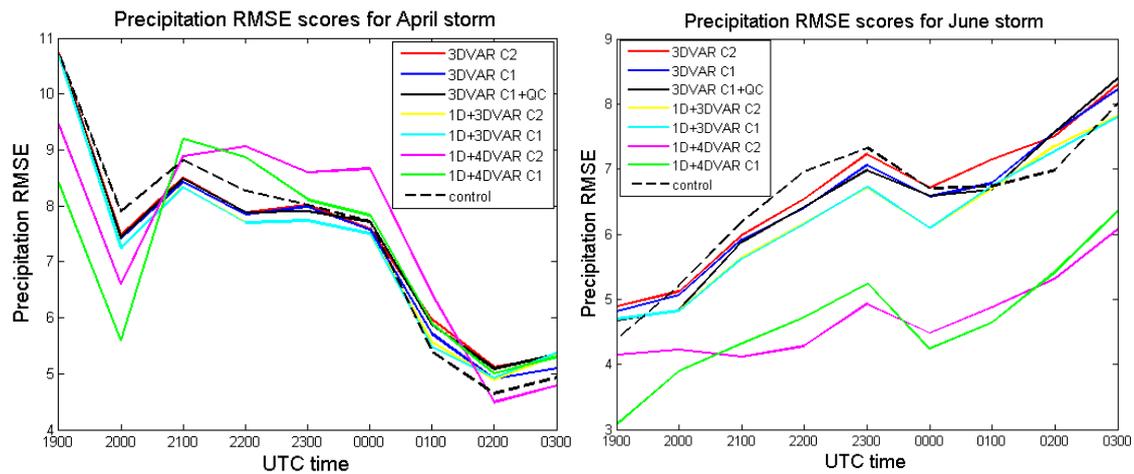

Figure 23. RMSE of precipitation (mm) for both study days compared with stage IV observations. Assimilation was not performed after 2000 UTC for the 1D+4D-VAR simulation, and not after 0000 UTC for the 3D-VAR and 1D+3D-VAR approaches.